\documentstyle[12pt]{article}
\newtheorem{definition}{Definition}
\newtheorem{theorem}{Theorem}

\newtheorem{lemma}{Lemma}

\newenvironment{proof}{\noindent {\em{Proof:}}}{\hfill $\Box$
\vspace{4mm}}

\newcommand{\la}{\lambda}
\newcommand{\al}{\alpha}
\newcommand{\be}{\beta}
\newcommand{\Om}{\Omega}
\newcommand{\si}{\sigma}
\newcommand{\ep}{\varepsilon}
\newcommand{\comb}[2]{\pmatrix{#1\cr #2}}
\newcommand{\om}{\omega}
\newcommand{\Ga}{\Gamma}
\newcommand{\Lam}{\Lambda}

\newcommand{\Spin}{\mbox{\rm Spin}}

\newcommand{\pul}{{{\frac{1}{2}}}}

\begin{document}
\title{Eigenvalues of conformally invariant operators on
spheres.\thanks{Research supported by the grant GA\v CR 201/96/0310}}
\author{Jarol\'{\i}m Bure\v{s} , Vladim\'{\i}r Sou\v{c}ek}
\date{}
\maketitle
\section{Introduction}

Let $M$ be a smooth oriented compact n-dimensional
manifold, $n\geq 3$,
endowed with a Riemannian metric and a spin structure.
A huge amount of information has been collected concerning
spectral properties of the basic invariant differential operators
on $M.$ Spectra of the Laplace and the Dirac operators
has been computed explicitely on many examples of homogeneous spaces
(\cite{B1,B2,BGM,CaH,CFG,Mi,MS,Sa}).
Estimates concerning the first eigenvalues and relation to
geometry of $M$ has been studied in many papers
(\cite{Fri1,Hi,Kir,KSW,Li2,Lo,Su})
In a general case, exact formulae for eigenvalues are not
available but their asymptotic behaviour is a classical
subject studied for a long time already (\cite{DF,Ga}).

Recently, a growing interest is paid to properties of more
complicated  invariant first order differential operators on
$M.$ A prototype of them is the Rarita-Schwinger
operator (see
\cite{Fra1,Fra2,FraS,MP,N,NGRN,Pe1,Pe2,Pe3,Pe4,Pe5,RaS,US,Wa}).
It acts on sections of the bundle associated
to a more complicated representation of the group $\Spin(n).$
In the paper, we are going  to study spectral properties
of a certain class of differential
operators on $M$ which has been intensively used in
Clifford analysis in connection with monogenic differential
forms (see \cite{DSS,Ry1,So1,So2,SoS,So3}). The aim of the  paper is to compute
explicitely spectra of this class of conformally invariant
operators on the flat model, i.e. on spheres.

As for the Dirac and the Laplace operators, methods of
representation theory can be used in homogeneous case.
The main tool used in the paper are general results
of Branson, \'Olafsson and \O rsted (see \cite{BOO})
describing a construction of intertwinning operators
between principal series representations of semisimple Lie
groups.
They are able to compute spectra of a wide class of invariant
operators up to a normalisation, i.e. they are giving explicit
formulae for ratios of eigenvalues. These formulae can be used
directly in odd dimensions. In even dimensions,
differential operators studied here are not covered by the
results in \cite{BOO}, nevertheless the methods used there
can be adapted for our purpose (see Sect.3).

The symbol of the Dirac operator is given by the Clifford
multiplication. Hence the question of normalisation is answered
here by a choice of the Clifford action.
For higher spin representations and the associated invariant
operators, the question of normalisation of the studied operators
is first to be settled (see Sect.2).
To compute exact formulae for spectra, it is then sufficient
to find explicitely one eigenvalue. It is done in Sect.4
using methods developed in \cite{VSe}.

\section{First order conformally invariant operators}
A classification and a description of
first order conformally invariant differential
operators was first described
by Fegan in (\cite{F}). There is a standard definition of an
invariant (homogeneous) operator on  homogeneous spaces but
there are several different definitions of
conformally   invariant operators in a curved case
(for details see a
\cite{BE,S,CSS1,CSS2}). A construction of curved analogues of
invariant operators is a difficult task which is not yet
completely understood (see \cite{BE,GJM}).
For first order operators, however, there are no additional complications
in the curved case with respect to the homogenous model.
A general scheme for a construction of such invariant operators
 is as follows (see \cite{F}).

Let M be a compact oriented  manifold with a conformal structure.
Let us choose a Riemannian metric in the given conformal class
and suppose that a spin structure is given on $M,$ i.e. that we
have
  principal
fibre bundles
$$
\tilde{\cal P}\equiv\tilde{\cal P}_{Spin}\rightarrow {\cal P}_{SO} \rightarrow M.
$$
on the manifold $M.$

Finite-dimensional irreducible representations ${\bf V}_{\la}$ of the group
$H=\Spin(n)$ are classified by their highest weights $\la \in
\Lambda^{+},$
                    where for $n=2k$ even, we have
$$
\Lambda^+ =
\{ \la =(\la_1,...,\la_k); \la_1\geq \la_2\geq ... \geq \la_{k-1}\geq
|\la_k|\},
\la_i \in {\bf Z}\cup \pul {\bf Z}
$$
 and for $n=2k+1$ odd, we have
$$
\Lambda^+ =
\{ \la =(\la_1,...,\la_k); \la_1\geq \la_2\geq ... \geq \la_{k-1}\geq
\la_k\geq 0\},
\la_i \in {\bf Z}\cup \pul {\bf Z}.
$$

Invariant operators are acting among spaces of sections of
the corresponding associated bundles
$$
V_{\la} = {\tilde{\cal P}}\times_{H} {\bf V}_{\la}
$$
over $M.$
Let us consider the Levi-Civita connection $\omega$ of the chosen
Riemannian metric on ${\cal P}$ and let $\tilde{\omega}$ be its (unique) lift
to ${\tilde{\cal P}}.$
For any choice of ${\la}\in\Lam^+,$ we have the associated covariant derivative
$$
\nabla_{\la}:\Gamma(V_{\la})
{\rightarrow} \Gamma(V_{\la}
\otimes T^*(M)).
$$
There are standard algorithms
(see \cite{Sal}) for a decomposition of the tensor
product ${\bf V_{\la}}\otimes {\bf C}_n$ into irreducible components
$$
{\bf V_{\la}}\otimes {\bf C}_n = \oplus_{\la'\in A} {\bf V}_{\la'},
$$
where $A$ is the set of  highest weights of all irreducible
components (multiplicities included).
There are simple rules how to describe $A=A(\la)$ explicitely for any $\la$
(see \cite{F,S}).
Let $\pi_{\la'}$ be the projection from
${\bf V}_{\la}\otimes {\bf C}_n$ to
${\bf V}_{\la'}.$ Then  operators
$$
D_{\la,\la'} :\Gamma(V_{\la}) \rightarrow \Gamma (V_{\la'}),\;
D_{\la,\la'}:=\pi_{\la'}\circ\nabla^{\la}
$$
are first order conformally invariant differential operators
and all such operators
can be  constructed in this  way.

Any conformally invariant first order differential
 operator is uniquely determined
(up to a constant multiple) by a choice of allowed $\la$ and $\la'$
but there is no natural normalization in general.
 To study spectral properties, it is necessary to remove this
 ambiguity and to fix a scale of the operator, to choose appropriate
 normalization. For the Dirac operator, the choice of normalization
 is given by the Clifford action. By using twisted Dirac operators,
 we shall extend this normalization to a wide class of first order
 operators.

\begin{definition}
 Let ${\bf S}^{}$ (for $n=2k+1$), resp.
 ${\bf S}^{}={\bf S}^{+}\oplus{\bf S}^{-}$
(for $n=2k$), denote the basic spinor representations
 with highest weights
$\si=({\frac{1}{2}},\ldots,{\frac{1}{2}},{\frac{1}{2}}),$
resp.
$\si^{\pm}=({\frac{1}{2}},\ldots,{\frac{1}{2}},\pm{\frac{1}{2}}).$

Let $\la
\in\Lam^+,$ (for $n=2k+1$), resp. $\la^{\pm}\in\Lam^+$ (for
$n=2k$)
 be  dominant weights with
$\la =(\la_1,...,\la_{k-1},{\frac{1}{2}}),$
resp. $\la^{\pm}=(\la_1,...,\la_{k-1},\pm{\frac{1}{2}}).$
 Denote further $\la'=\la-\si\in\Lam^+,$ resp.
${\la'}=\la^{+}-\si^{+}\in\Lam^+.$
In even dimensions, we shall use the notation
$$
{\bf V}_{\la}={\bf V}_{\la^+}\oplus {\bf V}_{\la^-}.
$$

The representation ${\bf V}_{\la}$  appears with multiplicity one
in the decomposition of the tensor product
${\bf S}\otimes {\bf V}_{\la'}$ (it is the Cartan product of both
representations). Hence we can write the product
as
$$
{\bf S}\otimes {\bf V}_{\la'}=
{\bf V}_{\la} \oplus {\bf W},
$$
where ${\bf W}$ is the sum of all other  irreducible components
in the decomposition.

Let  $D^T_{\la'}$ be the twisted Dirac operator
on $S\otimes V_{\la'}.$  If we write the operator $D^T_{\la'}$
in the block form as

\vskip 1.7cm
{\center{
\unitlength 1.00mm
\linethickness{0.4pt}
\begin{picture}(73,50)

\put(25,57){$\Gamma(S\otimes V_{\lambda'})$}
\put(60,57){$\Gamma(S\otimes V_{\lambda'})$}

\put(30,42){$\Gamma(V_{\lambda})$}
\put(65,42){$\Gamma(V_{\lambda})$}

\put(30,27){$\Gamma(W)$}
\put(65,27){$\Gamma(W)$}

\put(44,58){\vector(1,0){13.5}}
\put(42,43){\vector(1,0){17.5}}
\put(42,28){\vector(1,0){17.5}}
\put(42,31.3){\vector(2,1){17}}
\put(42,39.7){\vector(2,-1){17}}

\put(32,50){$\|$}
\put(67,50){$\|$}

\put(32,35){$\oplus$}
\put(67,35){$\oplus$}

\put(49,62){\makebox(0,0)[cc]{$D^T_{\lambda'}$}}
\put(49,47){\makebox(0,0)[cc]{$D_{\lambda}$}}

\put(16,52){\makebox(0,0)[cc]{$$}}
\put(40,52){\makebox(0,0)[cc]{$$}}

\end{picture}
}}

\vskip -1.7cm
\noindent
we have defined 4 new invariant operators, one of them being the operator
$$
D_{\la}:\Ga(V_{\la})\rightarrow\Ga(V_{\la}).
$$
Operators $D_{\la}$ defined in such a way will be called
{\it higher spin Dirac operators}.
\end{definition}

A certain subclass of invariant operators
discussed above have appeared often in discussions
of higher dimensional generalizations of holomorphic differential
forms (see \cite{DSS,So2}). They are arising in the following way.
Let us consider spinor valued differential forms, they are
coming as elements of the twisted de Rham sequence,

$$  \Gamma(S^{\pm})\stackrel{\nabla^S}{\rightarrow}
  \ldots
  \Gamma({\Omega}^{k}_c \otimes S^{\pm})\stackrel{\nabla^S}{\rightarrow}
  \ldots
\stackrel{\nabla^S}{\rightarrow}
  \Gamma({\Omega}^{n}_c \otimes S^{\pm})
$$
where $\nabla^S$ denotes the associated covariant  derivative
on spinor bundles extended to $S$-valued forms (see
\cite{So2,VSe}).

Every representation
$\Lambda^k({\bf C}_n)\otimes {\bf S}^{}$  can be split
into irreducible pieces. There are no multiplicities in the
decomposition, so the irreducible pieces are well defined.
For $k$ forms ($k\leq[n/2]$), there are  $k$ pieces in the
decomposition and the decomposition is symmetric with
respect to the action of the Hodge star operator.
The space of spinor valued $k$-forms
$ \Gamma({\Omega}^{k}_c \otimes S^{\pm})$
($k\leq [n/2]$) can be written
as the sum $\oplus_{j=1}^kE^{k,j}$
and it can be checked (see \cite{DSS,VSe,So2}) that
$E^{k,j}$ is the bundle associated with the representation
with the highest weight
$\la_j=({\frac{3}{2}},\ldots,{\frac{3}{2}},{\frac{1}{2}},\ldots,
{\frac{1}{2}},
\pm{\frac{1}{2}}),$ where the number $j$ indicates that
the component ${\frac{3}{2}}$ appears with multiplicity equal to $j.$
Signs $\pm$ at the last components are relevant only  in even
dimensions (more details can be found in \cite{VSe}).
The whole splitting can be described by the following triangle
shaped diagram (in odd dimensions, there are two columns of the
same length in the middle).

\begin{displaymath}
\begin{array}{ccccccccccccc}
  E^{0,0} &\stackrel{D_0}{\longrightarrow}&
  E^{1,0} &\stackrel{D_0}{\longrightarrow}&
  \ldots &     \stackrel{D_0}{\longrightarrow}&
  E^{k,0} &\stackrel{D_0}{\longrightarrow}&
  \ldots & \stackrel{D_0}{\longrightarrow}&
  E^{2k-1,0} &\stackrel{D_0}{\longrightarrow}&
  E^{2k,0} \\

  && \oplus && \oplus && \oplus && \oplus && \oplus &&  \\

  &&
  E^{1,1} &\stackrel{D_1}{\longrightarrow}&
  \ldots &          \stackrel{D_1}{\longrightarrow}&
  E^{k,1} &     \stackrel{D_1}{\longrightarrow}&
  \ldots & \stackrel{D_1}{\longrightarrow}&
  E^{2k-1,1}
  &&
   \\
   &&  &&   \oplus && \oplus && \oplus &&  && \\
   &&  && \ldots &\stackrel{D_{j}}{\longrightarrow}& \vdots
   &\stackrel{D_{j}}{\longrightarrow}& \ldots &&  && \\
    &&  &&  && \oplus &&  &&  && \\
    &&  &&  && E^{k,k} &&  &&  && \\
\end{array}
\end{displaymath}

The general construction of invariant operators described above
can be used
in the special case of spinor valued forms.
The  covariant
derivative $\nabla ^S$ restricted to $E^{k,j}$ and projected
to $E^{k+1,j'}$ is an example of this general construction.
It can be shown that if $|j-j'|>1,$ then the corresponding
invariant operator is trivial. We shall be mainly interested
in 'horizontal arrows', i.e. in operators $D_j$ given
by restriction to $E^{k,j}$ and projection to $E^{k+1,j}.$
They are indicated in the above scheme. The simplest cases among
them are well known. The operator $D_0$ is (a multiple of) the Dirac
operator. The operator $D_1$ is (an elliptic version of)
the operator called Rarita-Schwinger operator by physicists
(see \cite{EP,RaS,Wa}). All of them are elliptic operators  (see
\cite{So1}).
Note that all operators $D_j$ on the same row in the scheme above
cannot be identified without further comments.
 To compare them, it is necessary first
to choose an equivariant isomorphism among corresponding bundles.
Then
 they coincide  up to a constant multiple.

To compare the operators $D_j$ in the above
scheme with the higher spin Dirac operators (see Def.1),
we shall  choose a certain identification of the
corresponding source and target bundles.
 We shall do it for the first operator $D_j$ in the row.

Let us characterize an algebraic operator
$Y:\Gamma(\Omega^{k+1}_c\otimes S)\rightarrow
\Gamma(\Omega^{k}_c\otimes S)$
by a local formula
$$
Y(\om\otimes s)=-\sum_i\iota(e_i)\om\otimes e_i\cdot s,
$$
where $\{e_i\}$ is a (local)  orthonormal basis of
$TM$ and $\iota$ denotes the contraction of a differential form
by a vector. As shown in \cite{VSe}, the map
$Y:E^{k+1,j}\rightarrow E^{k,j}, j<k<[n/2]$ is an isomorphism.

The twisted Dirac operator $D^T$ maps the space $\Gamma(\Omega^{k}_c\otimes S)$
to itself.
In \cite{VSe}, it was proved that we have a relation
$\nabla\circ Y+Y\circ\nabla=-D^T.$
Let us denote the projection from
$\Omega^{k}_c\otimes S$ onto $E^{k,j}$  by $\pi_{k,j}.$
 Symbols $\tilde{D}_j,\,0\leq j<[n/2]$ will denote
operators
$$
\tilde{D}_j=Y\circ D_j=
\pi_{j,j}\circ Y\circ \nabla^S|_{E^{j,j}},
$$
 mapping the space of sections of $E^{j,j}$ to itself.
Then $Y|_{E^{j,j}}=0$ implies that
$$
\tilde{D}_j=
\pi_{j,j}\circ Y\circ\nabla^S|_{E^{j,j}}=
-\pi_{j,j}\circ D^T|_{E^{j,j}}=-D_{\la_j},
$$
where $D_{\la_j}$ is the higher spin Dirac operator corresponding
to the bundle $V_{\la_j},\,
\la_j=({\frac{3}{2}},\ldots,{\frac{3}{2}},{\frac{1}{2}},\ldots,
{\frac{1}{2}})$ (component ${\frac{3}{2}}$ appearing $j$
times).
More precisely, there are no signs in odd dimensions, while in
even dimension, $V_{\la_j}=V_{\la^+_j}\oplus V_{\la^-_j}.$
To compute spectrum of the higher spin Dirac operators,
it is hence sufficient to do it for $\tilde{D}_j.$

Now, we shall restrict our study to operators $\tilde{D}_j$
and we shall consider them on  spheres.
We would like to compute their spectra.
The spectrum of the Dirac
operator is well-known (see \cite{B2}).
\begin{lemma}
The eigenvalues of the Dirac operator on the sphere $S_n$ with
standard metric are

$$
\mu_l=\pm\left({\frac{n}{2}}+l\right);\;l=0,1,2,\ldots.
$$
with multiplicity
$$
2^{[\frac{n}{2}]} \comb{l+n-1}{l} .
$$
\end{lemma}

The main  result of the paper is given in the following theorem.

\begin{theorem}
  Let $D_{\la_j}=-\tilde{D}_j,\,0<j< n/2,$
 be the higher spin Dirac operators defined above,
considered on the sphere $S_n$ with the standard metric.
Then their eigenvalues are :

$$
\mu^1_l=\pm\left({\frac{n}{2}}+l\right);\;l=1,2,\ldots.
$$
with multiplicity

$$
2^{[\frac{n}{2}]}\comb{n+1}{j+1}\comb{l+n}{l-1}
\frac{(n-2j)(j+1)}{(l+j)(l+n-j)}
$$
and
$$
\mu^2_l=\pm\left[{\frac{n-2j}{n-2j+2}}\left({\frac{n}{2}}+l
\right)\right];\;l=1,2,\ldots.
$$
with multiplicity
$$
2^{[\frac{n}{2}]}\comb{n+1}{j}\comb{l+n}{l-1}
\frac{(n-2j+2)j}{(l+j-1)(l+n-j+1)}.
$$

\end{theorem}

The rest of the paper will be devoted to the proof of the
theorem.

\section{Ratios of eigenvalues}
A main tool for computation of eigenvalues will be taken
from the paper of Branson, \'Olafsson and \O rsted (see
\cite{BOO}). Their paper is designed to {\it construct}
invariant operators on homogeneous spaces in a diagonal
form.
They have developed a powerful method of study of invariant
operators (not necessarily differential ones!) using
representation theoretical methods.
They are prescribing the so called spectral function,
giving eigenvalues (up to a multiple)
 of an operator in question on suitably defined
finite dimensional spaces of eigenfunctions.
It applies to a broad class of homogeneous spaces, which
includes the conformally invariant operators considered above in
case
 of odd dimensions.
But operators $D_{\la_j}$ in even dimensional case are
explicitely excluded from consideration in their paper.
Our task here is different so that we can {\it compute}
ratios of eigenvalues using their method also in even dimensions.

The first thing to note is that eigenspaces of our operators
can be easily identified and described using representation
theory.
Let us consider the $n$-dimensional sphere $S^{n}$ as a homogeneous space
$$
S^n = G/P = K/H
$$
where
$G = {\rm Spin}_o(n+1,1)$ ,  $K ={\rm Spin}(n+1)$ is a maximal
compact subgroup of $G,$ $H=\Spin(n)$  and
$P$ is a (noncompact) maximal
parabolic $P = M A N$ with $M\subset K$.
The invariant metric $g$ on $S^n$ is   constructed by left
translation of the Killing form
$\tilde{B} = B/2n$, then $S^n$ has
 constant sectional curvature
$K= 1$ with respect to this metric.
We shall need here the compact picture $S^n=K/H$ only.

Let $\la\in\Lam^+$ and let $V_{\la}$ be the corresponding
homogeneous bundle on the sphere.
The  group $K$ acts on the space of sections
$\Ga(V_{\la})$ by the left regular representation.
 The group $K$ is compact, hence the space of sections can be
 decomposed to
corresponding isotypic components, which are finite-dimensional.
The main case considered in \cite{BOO} is the multiplicity one case,
when these isotypic
components are irreducible. Then by Schur lemma, any invariant
operator (when restricted to these components and
 acting among identical bundles)
 is a multiple of identity. Then these
components are eigenspaces of the operator.
To compute ratios of eigenvalues, the authors use
a suitable combination of Casimir operators called the spectrum
generating operator.

                 There are explicit formulas how to find
highest weights of  isotypic components appearing in the
decomposition. They are given by the so called branching rules,
which were carefully studied in representation theory.
In the conformal case needed below, they are given as follows.

Let us agree first the following notation.
Let  $\la\in\Lam^+(\Spin(n))$
and  $\al\in\Lam^+(\Spin(n+1)).$
The symbol $\al\downarrow\la$ is defined by the following relations:

1) Let $n=2k.$

$$
\al\downarrow\la \iff \al_1\geq \la_1\geq\al_2\geq \la_2\geq \ldots\geq \al_{k}\geq
|\la_k|.
$$

2) Let $n=2k+1.$

$$
\al\downarrow\la\iff \al_1\geq \la_1\geq\al_2\geq \la_2\geq \ldots\geq \al_{k}\geq
 \la_k\geq|\al_{k+1}|.
$$

If we consider now the space of sections $W=\Gamma(V_{\la}),\,
\la\in\Lam^+(H)$ as a $K$~-~modul, then
all isotypic components $V_{\al},\al\in\Lam^+(K)$ have
multiplcity at most one and are nontrivial
iff $\al\downarrow\la.$
The sum $\oplus_{\al,\al\downarrow\la}W_{\al}$ is then dense in
$\Gamma(V_{\la}).$

Methods and results of \cite{BOO} can be used to show

\begin{theorem}
Let $D_{\la},\,\la\in\Lam^+,|\la_k|={\frac{1}{2}},$
 be a higher spin Dirac operator
(see Def.1) and let $\mu,\mu'$ be its two different eigenvalues,
having both the same sign.
Let $W,$ resp. $W'$ be the corresponding spaces of eigenvectors.

Then there exist
isotypic components $W_{\al}, W_{\al'}$
with highest weights
$\al,\al'\in\Lam^+(K)$ such that
$W\subset W_{\al},\,W'\subset W_{\al'},$
and
$$
{\frac{\mu}{\mu'}}=
\Pi_{a=1}^{[{\frac{n+1}{2}}]}
{\frac
{\Gamma({\frac{1}{2}}(n+3)-a+\al_a)
\Gamma({\frac{1}{2}}(n+1)-a+\al'_a)}
{\Gamma({\frac{1}{2}}(n+3)-a+\al'_a)
\Gamma({\frac{1}{2}}(n+1)-a+\al_a)}
}
$$
\end{theorem}

\begin{proof}
Suppose first that $n$ is odd, $n=2k+1.$
Then the space ${\bf V}_{\la}$ is irreducible,
all isotypic conponents of $\Gamma(V_{\la})$ have
multiplicity one and the formula above for ratios of their eigenvalues
was proved in \cite{BOO}.

So suppose next that $n=2k.$
In the paper \cite{BOO}, they have to exclude this case from their
construction.
The main reason was that in this case,
they had no control over certain compatibility conditions
needed  for it.
Nevertheless, if the aim is not to {\it construct} intertwining
operators but to {\it compute} their eigenvalues under the assumption
that they exist, it is not necessary to verify corresponding
compatibility conditions and their methods are applicable.

So only problem to discuss is that isotypic components of the space of
sections
$\Gamma(V_{\la})$ have multiplicity two.
Indeed,
${\bf V}_{\la}={\bf V}_{\la^+}\oplus {\bf V}_{\la^-}.$
Isotypic components of $\Gamma(V_{\la^+}),$ have all multiplicity
one, but they are identical with the corresponding isotypic components of
$\Gamma(V_{\la-}).$

The operator $D_{\la}$ intertwines action of $K,$ so it preserves
individual isotypic components.
If $s=(s^+,s^-),\,s^{\pm}\in\Gamma(V_{\la^{\pm}})$ is eigenvector
of $D_{\la}$ with eigenvalue $\mu,$ then $(s^+,-s^-)$ is
eigenvector with eigenvalue $-\mu.$Hence the restriction of $D_{\la}$ to any isotypic
component $W_{\al}$ has  at least two eigenvalues $\pm\mu.$
The corresponding eigenspaces are then $K$-modules
and the isotypic component $W_{\al}$ is a direct sum of them.
Let us denote by $W^+$ the closure of the sum of all eigenspaces
corresponding to positive eigenvalues.
The corresponding isotypic components have multiplicity one,
$W^+$ is an invariant subspace with respect to the action of $G$
and the computation in \cite{BOO} can be repeated to
prove the result.

\end{proof}

\section{Normalisation of eigenvalues}
To finish the computation of spectra, it is necessary to compute
at least one eigenvalue of a given operator.
The spectrum of the Dirac operator is known.
We shall show how to compute inductively one eigenvalue for
operators $\tilde{D}_j=-D_{\la_j},\,0<j<[n/2].$
It will lead then in next section to a formula for their full
spectrum.

A useful relation among spectra operators $\tilde{D}_j$
was shown in \cite{VSe}, the following
theorem is proved there.

\begin{theorem}
Let $M$ be an Einstein spin manifold.
Let us define a first order differential
operator $T_j$  by
$$T_j=\pi_{j+1,j+1}\circ\nabla^S|_{E^{j,j}},\,0\leq j\leq[n/2]-1.$$

If $s$ is an eigenvector of the operator $\tilde{D}_j$ corresponding
to an eigenvalue~$\mu$ and if $T_j(s)\not =0,$
then $s'=T_j(s)$ is an eigenvector of the operator $\tilde{D}_{j+1}$
corresponding to the eigenvalue $\mu'=-{\frac{n-2}{n}}\mu.$

\end{theorem}

As a consequence, if we are able to find eigenvectors of $\tilde{D}_j$
which does not belong to the kernel of $T_j,$ we can compute
at least one eigenvalue of $\tilde{D}_{j+1}.$
The following theorem shows that it is always possible.

\begin{theorem}
The operators $T_j,\,0\leq j\leq \frac{n}{2}-1$ have nontrivial symbol,
hence their kernel is a proper subset of $\Gamma(V_{\la_j}).$
\end{theorem}

\begin{proof}
 Let $\ep(v):\Om^j_c\rightarrow \Om^{j+1}_c,\,v\in \Om^1_c$
 denote the outer multiplication by the element $v.$
 Then symbol $\si$ of the operator $\tilde{D}^j$
 is given by
 $$
 \si(v)(\om)=\pi_{j+1,j+1}\circ\ep(v)(\om),\;
 v\in\Om^1_c,\,\om\in E^{j,j}\subset\Om^j\otimes S.
 $$

Let $v_j$ denote a nontrivial weight vector of the fundamental
representation ${\bf C}^n$ corresponding to a weight $\la_j=
(0,\ldots,1,\ldots, 0)$ with $1$ on the $j$-th place, resp.
corresponding associated element in $\Om^1_c.$
Then $v_1\wedge\ldots\wedge v_j$ is a (nontrivial) weight vector
of $\Lam^j{\bf C}^n.$ Denote further by $s_0$ a nontrivial weight vector
for the highest weight of ${\bf S}.$ Then
$$
w=v_1\wedge\ldots\wedge v_j\otimes s_0
$$
is a (nontrivial) weight vector for the highest weight of
$\Lam^j{\bf C}^n\otimes {\bf S}.$
 Hence $w$ belongs to the Cartan product of $\Om^j\otimes S,$
 which is just equal to $E^{j,j}.$

Hence
$$
\si(v_{j+1})(v_1\wedge\ldots\wedge v_j\otimes s_0)=
(-1)^jv_1\wedge\ldots\wedge v_{j+1}\otimes s_0
$$
is a nontrivial vector and the theorem is proved.
\end{proof}

\section{The proof of Theorem 1}
Now we can finish the proof of the main theorem.
It is necessary to distinguish even and odd dimensional cases.

\begin{proof}

\noindent 1) Let first $n=2k+1.$

The highest weight of the space $S=E^{0,0}$ is
$\la_0=({\frac{1}{2}},\ldots,{\frac{1}{2}})$ ($k$ components)  and the space
of sections of $\Gamma (E^{0,0})$ is a sum
of $K$-types
$$
A_{\al_0(\pm,l)},\;
\al_0(\pm,l)=
\left({\frac{2l+1}{2}},
{\frac{1}{2}},\ldots,{\frac{1}{2}},
\pm{\frac{1}{2}}\right),\,l=0,1,2,\ldots ,
$$
($\al_0(\pm,l)$ having $k+1$ components).

The highest weight of the space $E^{j,j},\,j>0$ is
$ \la_j=\left({\frac{3}{2}},\ldots,{\frac{3}{2}},
 {\frac{1}{2}},\ldots,{\frac{1}{2}}\right)$ (with ${\frac{3}{2}}$
 appearing $j$ times) and the space
of sections of $\Gamma (E^{j,j})$ is a sum
of $K$-types
$$
A_{\al_j(\pm,l)},
\;\al_j(\pm,l)=
\left({\frac{2l+1}{2}},
{\frac{3}{2}},\ldots,{\frac{3}{2}},
{\frac{1}{2}},\ldots,{\frac{1}{2}},\pm{\frac{1}{2}}\right),
\,l=1,2,\ldots ,
$$
and
$$
B_{\be_j(\pm,l)},
\;\be_j(\pm,l)=
\left({\frac{2l+1}{2}},
{\frac{3}{2}},\ldots,{\frac{3}{2}},
{\frac{1}{2}},\ldots,{\frac{1}{2}},\pm{\frac{1}{2}}\right),
\,l=1,2,\ldots ,
$$
where the component ${\frac{3}{2}}$ is appearing $j-1$ times
in the weight $\al_j(\pm,l)$ and $j$ times in the weight $\be_j(\pm,l).$

Using the formula for the ratio of eigenvalues
 from Th.2, we get first for $\al=\al_j(\pm,l)$
$$
\Pi_{a=1}^{{\frac{n+1}{2}}}
{\frac
{\Gamma({\frac{1}{2}}(n+3)-a+\al_a)}
{\Gamma({\frac{1}{2}}(n+1)-a+\al_a)}
}=
\pm({\frac{n}{2}}+l)2^{-{\frac{n-1}{2}}}(n!!)({\frac{n}{2}}-j)^{-1},
$$
hence
 there are constant $C_1,C_2$
 (independent of a $K$ type chosen, but depending on $n$ and $j$)
 such that the eigenvalues
 $\mu^1_{\pm,l}(j),$ resp.  $\mu^2_{\pm,l}(j),$
 corresponding to the eigenspace $A_{\al(\pm,l)},$
resp.   $B_{\be(\pm,l)}$
are equal to
$$
\mu^1_{\pm,l}=\pm C_1\,(l+{\frac{n}{2}})
$$
resp.
$$
\mu^2_{\pm,l}=\pm C_2\,(l+{\frac{n}{2}}).
$$
Moreover, Th.2 implies that
$$
{\frac{\mu^1_{\pm,l}(j)}{\mu^2_{\pm,l}(j)}}=
{\frac
{n-2j}{n-2j+2}
}.
$$

The unknown constants $C_1,C_2$ will be computed inductively (with respect
to~$j$). For the Dirac operator, the spectrum is known (see e.g.
\cite{B2}),
the eigenvalue corresponding to $K$-typ with $\al(\pm,l)$ is equal
to $({\frac{n}{2}}+l).$

For the Rarita-Schwinger operator ($j=1$), there are two sequences of
$K$-types, one of them being a subset of that for the Dirac operator
(only the first term is missing).

 The twistor operator $T_0$ is invariant (hence should preserve the
 label of a $K$-type)
and has a finite dimensional kernel (hence is nontrivial for at least
one $K$-type). Th.3 is then saying that
$$
{\mu^1_{\pm,l}(1)}= \pm\left({\frac{n}{2}}+l\right){\frac{n-2}{n}},
$$ hence the theorem is valid for $j=1.$

Due to the preceding theorem, operators $T_j$ are nontrivial for all
$j,$ so the proof can be finished in the same way by induction.

\vskip 2mm
\noindent 2) Let $n=2k.$
In even dimensions, $E^{j,j},j>0$ is a sum of two $\pm$ spaces with
highest weights
$\la_j^{\pm}=\left({\frac{3}{2}},\ldots,{\frac{3}{2}},
 {\frac{1}{2}},\ldots,{\frac{1}{2}},\pm{\frac{1}{2}}\right)$ ($k$ components,
  ${\frac{3}{2}}$
 appearing $j$ times).
 Hence the space of sections will be a sum of $K$-types
 ($\al$'s having also $k$ components)
$$
A_{\al(l)},\al(l)=
\left({\frac{2l+1}{2}},
{\frac{3}{2}},\ldots,{\frac{3}{2}},
{\frac{1}{2}},\ldots,{\frac{1}{2}}\right),
\,l=1,2,\ldots ,
$$
and
$$
B_{\be(l)},\be(l)=
\left({\frac{2l+1}{2}},
{\frac{3}{2}},\ldots,{\frac{3}{2}},
{\frac{1}{2}},\ldots,{\frac{1}{2}}\right),
\,l=1,2,\ldots ,
$$
where  ${\frac{3}{2}}$ is appearing $j-1$ times
in the weight $\al(l)$ and $j$ times in the weight $\be(l).$
This time, however, each type will appear with a multiplicity two.

As in the proof of Th.2, we can split each isotypic component
with respect to $K$ as a sum of eigenspaces corresponding to
opposite eigenvalues. The sum of spaces corresponding to
positive ones will be invariant with respect to $G$ and
the same proof as in odd dimensional case will go through.

The formula for the dimension of the space of eigenvectors
is the consequence of the Weyl dimensional formula for
the representation with the corresponding highest weight.
\end{proof}

\markright{References}

\end{document}